\documentclass[12pt,a4]{article}
\usepackage{amssymb,amsfonts,amsmath,setspace, color} 

\renewcommand\le{\leqslant}
\renewcommand\ge{\geqslant}

\newcounter{Remark}

\newtheorem{Theorem}{Theorem}

\newtheorem{Remark}{Remark}

\newcommand{\f}{\varphi}
\newcommand{\si}{\sigma}
\newcommand{\de}{\delta}
\newcommand{\D}{\Delta}                  

\newcommand{\exn}{{\bf E}\,}
\newcommand{\pr}{{\bf P}}
\newcommand{\R}{{\mathbb{R}}}
\newcommand{\C}{{\mathbb{C}}}
\newcommand{\Z}{{\mathbb{Z}}}

\begin{document}
%
%
%
%

\title{A refined version of the integro-local Stone theorem$^1$}

\author{Alexander A.\ Borovkov$^2$ and  Konstantin A.\ Borovkov$^3$}

\date{}
\maketitle

\footnotetext[1]{Research supported by the President of the Russian Federation Grant NSh-3695.2008.1,  the  Russian Foundation for Fundamental Research Grant 08--01--00962 and the ARC grant DP150102758.}

\footnotetext[2]{Sobolev Institute of Mathematics,  Ac.~Koptyug avenue~4, 630090 Novosibirsk, Russia. E-mail: borovkov@math.nsc.ru.}

\footnotetext[3]{School of Mathematics and Statistics, The University of Melbourne, Parkville 3010, Australia. E-mail: borovkov@unimelb.edu.au.}

\begin{abstract}
Let $X, X_1, X_2,\ldots $ be a sequence of non-lattice  i.i.d.\ random variables  with $\exn X=0,$ $\exn X=1,$ and let   $S_n:= X_1+ \cdots+ X_n$, $n\ge 1.$ We refine Stone's integro-local theorem by deriving the first term in the asymptotic expansion for the probability $\pr \bigl(S_n\in [x,x+\D)\bigr)$ as $n\to\infty$ and establishing uniform bounds for the remainder term, under the assumption that the distribution of $X$ satisfies Cram\'er's strong non-lattice condition and $\exn |X|^r<\infty$ for some $r\ge 3$.


\end{abstract}

\section{Introduction and the main result}

In the present note, we establish  a refinement of the following remarkable integro-local version of the central limit theorem due to C.~Stone \cite{St65, St67}. Let $X, X_1, X_2,\ldots $ be a sequence of non-lattice independent identically distributed (i.i.d.) random variables  (r.v.'s) following a common distribution~$F$ such that $\exn X=0,$ $\exn X^2=1,$ and let   $S_n:= X_1+ \cdots+ X_n$, $n\ge 1.$ For  $x\in \R$ and $\D>0$, set
\[
\D[x) := [x, x+\D).
\]
Then, as $n\to \infty$,
\begin{align}
\frac1\Delta \pr \bigl(S_n\in \D [x)\bigr)
   =    n^{-1/2} \phi (x  n^{-1/2})+ o(n^{-1/2}),
  \label{StoneThm}
\end{align}
where $\phi (t):=(2\pi)^{-1/2} e^{-t^2/2 }$  is the standard normal density and the remainder term is uniform in $x\in\R$ and in $\D\in [\D_0, \D_1]$ for any fixed $0<\D_0 < \D_1 <\infty$ (in fact, C.~Stone established  more general versions of the above result, including  convergence to stable laws, the multivariate case and large deviations).

It is quite appropriate  to call relations of the form \eqref{StoneThm} the {\em integro-local theorems}, to distinguish them from the integral theorems (which refer to approximating probabilities of the form $\pr (S_n < x)$, $x\in \R$) and the local ones (which deal  with approximating the densities of $S_n$ in the ``smooth" case; note that in the arithmetic case, the integro-local theorems are in fact the local ones: they concern  approximating probabilities $\pr (S_n=x)$ for $x\in \Z$).

The integro-local theorem is perhaps the most perfect and precise version of the classical central limit theorem. Indeed, it does not assume any additional conditions on top of the standard requirement of finite second moments (except for distinguishing between the lattice and non-lattice cases), but has basically got the same accuracy as the local limit theorems (note that, for small~$\D$, the left-hand side of~\eqref{StoneThm} is ``almost the density" of $S_n$),  without making any assumptions about existence of the densities.
That means that the  integro-local theorems are much more precise than the integral ones, and it is easy to see that one can derive the  assertions of the latter  from the former, but not the other way around.

The integro-local   theorems are rather effective and often the most adequate technical tools  in a number of problems in probability theory.   For instance, they are used  for computing the exact asymptotics of  large deviation probabilities for sums of independent r.v.'s (cf.\ Chapter~9 in~\cite{Bo13}). They also proved instrumental for studying the distribution of the first passage time of a curvilinear boundary by a random walk~\cite{Bo16}, establishing integro-local theorems for compound renewal processes (see Chapter~10 in~\cite{Bo13}) and in a number of other problems.

Concerning the history of the problem, note that a special case of relation~\eqref{StoneThm} (when $x$ is fixed) was first established  by L.A.~Shepp~\cite{Sh64}. A textbook exposition of the proof of~\eqref{StoneThm} can be found in Section~8.7 of~\cite{Bo13}. Under additional Cram\'er's conditions (the moment generating function is finite in a neighborhood of zero and the strong non-lattice condition on the characteristic function of~$X$ is met, see~\eqref{CC} below), relation~\eqref{StoneThm} was extended in~\cite{BoMo99, BoMo01} in the multivariate setting  to an asymptotic expansion in the powers of $n^{-1/2}$ and also to the case where $\D_0$ can be vanishing.  Extensions of Stone's theorem to the case of  non-identically distributed independent r.v.'s in the triangular array scheme (covering the large deviations zone as well) were established in~\cite{Bo11}.

In the lattice case, an analog of~\eqref{StoneThm} was obtained by B.V.~Gnedenko in the univariate case (see Chapter~9 in \cite{GnKo49}) and by E.L.~Rvacheva~\cite{Rv54}  in the multivariate case.

It is most natural to ask if the remainder term in~\eqref{StoneThm} can be sharpened under minimal additional assumptions. The first step in that direction was made in~\cite{Bo16a}, where it was shown that, under Cram\'er's strong non-lattice  condition
\begin{equation}
\limsup_{|\lambda|\to \infty} |\f (\lambda)|<1
\label{CC}
\end{equation}
on the  characteristic function (ch.f.)
$
\f (\lambda):=\exn e^{i\lambda X},$ $\lambda \in\R,
$
of $X$, and the moment condition $\exn |X|^r<\infty$ for some $r \in (2,3],$ relation~\eqref{StoneThm} holds with the reminder term replaced by $O(\D n^{-(r-1)/2})$ uniformly in $x\in\R$ and in $\D\in (q^n, cn^{(3-r)/2})$ for some $q\in (0,1)$ and every fixed $c>0$. In fact, \cite{Bo16a} actually establishes a multivariate version of that result.

In the present note, we further develop the approach from~\cite{Bo16a} to derive the first term of the asymptotic expansion for $\pr (S_n\in \D [x))$ with uniform bounds for the remainder term in the case when condition~\eqref{CC} is met and  $\exn |X|^r<\infty$ for some $r\in [3,4].$ It will be seen from the proofs that, under appropriate moment conditions, one can extend these results to asymptotic expansions with more terms. However, since  the very form of such expansions and their derivations are getting quite cumbersome, while technically they are not much different from the one-term case, we will restrict ourselves to presenting the latter only.

To formally state our main results, we will need some further  notations. For $r\in (2,\infty)$ and $b\in (1,\infty],$ introduce the class $\mathcal F_{r,b}$ of distributions $F$ on $\R$ satisfying the following moment conditions: for $X\sim F$, one has  $\exn X=0$, $\exn X^2 =1$ and
\[
\exn |X|^r <b   .
\]
In particular, $\mathcal F_{r,\infty}$ is the class of all  zero mean unit variance distributions with a finite $r$th absolute moment. For $F\in \mathcal F_{3,\infty}$,  we set
\[
\mu_3 := \exn X^3.
\]

Further, for  $\rho \in (0,1]$ and $b<\infty$, we denote by  $\mathcal F_{r,b}^\rho$ the totality of distributions from $\mathcal F_{r,b}$ that satisfy
\begin{equation}
\sup_{|\lambda|>1/b} |\f (\lambda)|<\rho.
\label{CC0}
\end{equation}
When $b=\infty$, we will understand by  $\mathcal F_{r,\infty}^1$ just the totality of distributions from $\mathcal F_{r,\infty}$ that satisfy
Cram\'er's strong non-lattice  condition~\eqref{CC}
(or, equivalently, $\sup_{|\lambda| >\varepsilon}|\f (\lambda)|<1 $ for any $\varepsilon >0$).


\begin{Theorem}
\label{Thm1}
{\rm (i)}~For any distribution $F\in\mathcal F_{3,\infty}^1$,  one has
\begin{align}
\frac1\Delta \pr \bigl(S_n\in \D [x)\bigr)
 =\frac{1}{n^{1/2}}\,\phi \Bigl(\frac{x}{n^{1/2}}\Bigr)
 \biggl(1+\frac{\mu_3 x}{6n}\biggl(\frac{x^2}{n}-3\biggr)- \frac{\D x}{2n}\biggr)
 +  \frac{R_n}n  ,
 \label{asym_form}
\end{align}
where for the remainder term $R_n=R_n(x, \D)$ the following holds true: there exists a $q\in (0,1)$ such that, for any fixed $\D_1>0$,
\[
\lim_{n\to\infty} \sup_{q^n \le \D \le \D_1}\sup_{x\in \R}  |R_n(x, \D)|=0.
\]
{\rm (ii)}~Moreover, for any fixed $r\in (3,4],$ $b<\infty$ and $\rho \in (0,1),$ representation \eqref{asym_form} holds true with the following uniform remainder bound over the  distribution class $ {\mathcal F}_{r,  b}^\rho:$  there exists a $q\in (0,1)$ such that, for any fixed $\D_1>0$,
\[
\limsup_{n\to\infty} \sup_{F \in {\mathcal F}_{r, b}^\rho} \sup_{q^n \le \D \le \D_1} \sup_{x\in \R}  n^{(r-3)/2} |R_n(x, \D)|<\infty .
\]
\end{Theorem}

\begin{Remark}
\rm  Note that the lower bound $q^n$ for the range of $\Delta$ values in the theorem cannot be ``qualitatively" improved. Indeed, assume that $F$ satisfies the conditions of part~(i) and has an atom of size $q_0\in (0,1)$ at zero. Then $\pr (S_n=0)\ge q_0^n$. Therefore, for $\Delta\in (0,q_0^n)$, the left-hand side of~\eqref{asym_form} will be at least one, whereas the right-hand side of that relation will be~$o(1),$ so that~\eqref{asym_form} cannot hold true for such values of~$\Delta$.
\end{Remark}

\begin{Remark}
\rm  In the general case of non-standardized r.v.'s $X$ with some $\mu:=\exn X $ and $\si^2 >0,$ expansion~\eqref{asym_form} will hold with $x$, $\D$ and $\mu_3$ on its right-hand side  replaced by $\si^{-1} (x-\mu),$ $\si^{-1}\D $ and $\si^{-3}\exn (X-\mu)^3$, respectively. The uniformity assertion in part~(ii) will have to be reformulated in that case as well.
\end{Remark}

\begin{Remark}
\rm
Note that our Theorem~\ref{Thm1}  can be viewed as an integro-local version  of the famous Chebyshev--Cram\'er asymptotic expansion   (a.k.a.\ the Edgeworth expansion, due to the contribution made in~\cite{Ed05}) for the distribution function of $S_n$, which was introduced in~\cite{Ch87} and formally proved in~\cite{Cr28, Es45}. Assume for a moment that $\exn X^4<\infty$ and suppose that condition~\eqref{CC} is met.
Then, denoting by $\Phi$ the standard normal distribution function, by ${\rm He}_k (x):=e^{-D^2/2}x^k$, where $D:=\frac{d}{dx}$, the $k$-th Chebyshev--Hermite polynomial, and by $\gamma_k$ the $k$-th cumulant  of $X$,  $k=1,2,\ldots,$  one has the following asymptotic expansion (see e.g.\ Section~5.7 in~\cite{Pe95}): as $n\to\infty$,
\[
\pr \biggl(\frac{S_n}{n^{1/2}} <v \bigg) = \Phi (v) - \phi (v)\biggl[
 \frac{\gamma_3}{6 n^{1/2}} {\rm He}_2 (v)
 +\frac1n \biggl(
 \frac{\gamma_3^2}{72}  {\rm He}_5 (v) +\frac{\gamma_4}{24}  {\rm He}_3 (v)
 \biggr)
 \biggr] + o (n^{-1})
\]
uniformly in $v\in\R$. Taking a fixed $\D>0$, setting $x:=v n^{1/2}$ and expanding
the  terms in the expression  on the right-hand side with  $v$ substituted by $v+\D n^{-1/2}$, we obtain
\begin{align*}
\pr \bigl(S_n\in \D [x)\bigr)
 & =
 \pr \biggl(\frac{S_n}{n^{1/2}} <v +\frac{\D}{n^{1/2}}\bigg) -\pr \biggl(\frac{S_n}{n^{1/2}} <v \bigg)
 \\
 & =  \frac{\D}{n^{1/2}}\, \phi (v)
 \biggl(1+\frac{\gamma_3  }{6n^{ 1/2} }  {\rm He_3}(v)  - \frac{\D }{2n^{ 1/2}}{\rm He}_1 (v)\biggr) + o (n^{- 1})
 \\
 & =\frac{\D}{n^{1/2}}\, \phi (v)
 \biggl(1+\frac{\mu_3  }{6n^{ 1/2} } (v^3-3v)  - \frac{\D v}{2n^{ 1/2}}\biggr) + o (n^{- 1}),
\end{align*}
where we replaced $\gamma_3$ with $\mu_3$ as these two quantities coincide for standardized r.v.'s.  In this note, we prove that  the above relation holds already when $\mu_3$ is finite. Moreover,  we establish   that the  above asymptotic relation  holds uniformly in $\D\in [q^n, \D_1]$  for some $q\in (0,1)$ and any fixed $\D_1 >0$. That result   cannot be derived from any asymptotic expansions for distribution functions. In addition, we  give uniform bounds for the remainder term in the case where $\exn |X|^r <b<\infty $ for some $r\in (3,4]$.
\end{Remark}

\begin{Remark}
\rm  One might also mention here the nonuniform error bounds for asymptotic expansions for probabilities $\pr (n^{-1/2}S_n\in A)$ for convex sets~$A$ obtained in the multivariate case in~\cite{Fo82}. However, for the purposes of solving the problem addressed in the present paper, these bounds are not better than the standard approximation rates from the classical asymptotic expansions mentioned in the previous remark.
\end{Remark}

\section{The proof of the main result}

The proof uses the method of characteristic functions with smoothing. We will only prove part~(ii), since  the proof of part~(i) follows exactly the same scheme and is somewhat simpler as it requires fewer explicit bounds.  So we will assume throughout this section that  $F\in  {\mathcal F}_{r,  b}^\rho $ for some fixed $\rho\in (0,1),$ $r\in (3,4]$ and $b<\infty$.

First we will derive fine asymptotics for the ``smoothed" distributions of $S_n$, namely, for the distributions of
\[
\widetilde{S}_n := S_n - \de U, \quad  n\ge 1,
\]
where  $\de=\de (\D, n)  $ will be chosen later and  $U $ is an r.v.\ uniformly distributed over $(0,1)$ and independent of $\{X_n\}_{n\ge 1}$.
Denote by
\[
\psi (\lambda):= \exn e^{-i\lambda U} = \frac{1-e^{-i\lambda}}{i\lambda}, \quad \lambda \in\R,
\]
the ch.f.\ of $-U$. Since the ch.f.\ of   $\widetilde{S}_n $ is clearly  equal to  $\f^n (\lambda) \psi (\de \lambda)$ and the function
\[
|\f^n (\lambda) \psi (\de \lambda) \psi (\D \lambda)|\le \min\{1,4 \de^{-1}\D^{-1}\lambda^{-2}\}, \qquad \lambda \in \R,
\]
is integrable on $\R$, we have from the standard inversion formula for ch.f.'s that yields the increments of the respective distribution functions (see e.g.\ (3.11) in Section~XV.4 of~\cite{Fe67}) that
\begin{align}
 \pr \bigl(\widetilde{S}_n \in \D [x)\bigr)
   & =\frac{\D}{2\pi }\int  e^{-i\lambda x} \f^n (\lambda) \psi (\de \lambda) \psi (\D \lambda) d\lambda
   =  \frac{\D}{2\pi }(I_1 + I_2 + I_3),
\label{iii}
\end{align}
where we defined $I_j:= \int_{B_j} e^{-i\lambda x} \f^n (\lambda) \psi (\de \lambda) \psi (\D \lambda) d\lambda,$ $j=1,2,3,$ as the integrals over the respective regions
\begin{align*}
 B_1  &:= \{\lambda\in \R: |\lambda|<h_1 n^{-1/3}\}, 
 \\
 B_2 &:= \{\lambda\in \R: h_1  n^{-1/3}\le |\lambda|< h_2 \},
 \\
 B_3 &:= \{\lambda\in \R: |\lambda|\ge  h_2  \}
\end{align*}
for  fixed $h_j >0,$ $j=1,2$,  to be chosen later.

Letting
\[
u:=\lambda n^{1/2},\quad v:= x n^{-1/2},\quad \widehat{\psi } (\lambda):= \psi (\de \lambda) \psi (\D \lambda),
\]
we can re-write $I_1$ as
\begin{align}
I_1= n^{-1/2}\int_{A_n}  e^{-i u v } \f^n (un^{-1/2})
\widehat{\psi} (un^{-1/2}) d u, \quad A_n:=(-h_1 n^{1/6}, h_1 n^{1/6}).
\label{I_1}
\end{align}
Our first step in evaluating $I_1$ consists in deriving representation~\eqref{f^n} below  for the second factor in the integrand. Make use of following expansion for the ch.f.\ $\f$:
\begin{align}
1-\f (\lambda) = \frac{1}2 \lambda^2 + \frac{i\mu_3}6\lambda^3 + \theta (\lambda) \lambda^3,
\label{theta}
\end{align}
where
\begin{align}
 |\theta (\lambda)|\le\frac{2^{4-r}\exn|X|^r|\lambda|^{r-3}}{r(r-1)(r-2)}
  \le \frac{b }{3} |\lambda|^{r-3}
\label{theta_b}
\end{align}
(see e.g.\ Section~12.4 in~\cite{Lo55}). Note that, by Lyapunov's inequality, for $|\lambda|\le b^{-1}$ one has
\begin{align}
\biggl|\frac{i\mu_3}6\lambda^3 + \theta (\lambda) \lambda^3\bigr|
 \le \frac{b^{3/r}}{6} |\lambda|^3 + \frac{b}{3}|\lambda|^{ r}
  = \frac{|b^{1/r} \lambda|^3}{6}   + \frac{|b^{1/r}\lambda|^{ r}}{3}
   \le \frac{b^{3/r-1}}2 \lambda^2.
   \label{la^2}
\end{align}
Hence, for $\lambda$ from that range,
\begin{align}
|1-\f (\lambda)| \le \frac{\lambda^2}2 + \frac{b^{3/r-1}}2 \lambda^2
 \le \lambda^2 \le b^{-2}<1.
 \label{1-f}
\end{align}
Since by the Taylor formula with remainder in  Lagrange form one has
\begin{align}
|\ln (1- z) + z|\le \frac{|z|^2}{2(1-c)^2} , \quad z\in\C, \ |z|\le c < 1,
 \label{ln (1-z)}
\end{align}
we conclude  that,   in the domain $|u|< b^{-1} n^{1/2}$,
\begin{align}
n\ln \f (un^{-1/2}) & =n\ln \bigl( 1 - (1-\f (un^{-1/2})\bigr)
 \notag \\
& = n\ln \biggl( 1 - \frac{u^2}{2n}  - \frac{i\mu_3u^3}{6n^{3/2}}
 +   \frac{u^3}{ n^{3/2}}\theta(un^{-1/2})\biggr)
 =  - \frac{u^2}{2 }+w,
\label{ln_f}
\end{align}
where
$
 w: =      - \frac{i\mu_3u^3}{6n^{1/2}}
 +   \frac{u^3}{ n^{1/2}}\theta_1(un^{-1/2})
$
and, in view of~\eqref{theta_b}, \eqref{1-f} and \eqref{ln (1-z)}, for $ |\lambda |\le b^{-1}$ one has
\begin{align}
|\theta_1 (\lambda) | \le |\theta  (\lambda) |
+ \frac{|1-\f(\lambda) |^2}{2 (1-b^{-2})^2|\lambda|^3}
 \le \frac{b}{3} |\lambda |^{r-3} + \frac{| \lambda  | }{2 (1-b^{-2})^2 }
 .
\label{theta_1}
\end{align}

Now from~\eqref{ln_f} we have
\begin{align*}
\f^n (un^{-1/2})  & =e^{-u^2/2}e^{w}
 =e^{-u^2/2}\bigl[1 + w + (e^{w} -(1+w))\bigr]  ,
\end{align*}
where,  again from the Taylor formula with remainder in  Lagrange form,
\[
|e^{w} -(1+w)|\le \frac{|w|^2}{2}e^{|w|}.
\]
Setting $h_1:= \bigl(
\frac{b^{3/r}}6 +\frac{b}3 + \frac1{2(1-b^{-2})^2}\bigr)^{-1/3} $, we have from~\eqref{theta_1} that, for $u\in A_n$,
\[
|w| =\biggl|
-  \frac{i\mu_3u^3}{6n^{1/2}}
 +   \frac{u^3}{ n^{1/2}}\theta_1(un^{-1/2})
\biggr|\le h_1^{-3}\frac{|u|^3}{n^{1/2}} <  1.
\]
From here and the previous two displayed formulae we see that, for $u\in A_n$ with the chosen $h_1$,  one has
\begin{align}
\f^n (un^{-1/2})  & =   e^{-u^2/2} \biggl( 1  - \frac{i\mu_3u^3}{6n^{1/2}}
 +   \frac{u^3}{ n^{1/2}}\theta_2(un^{-1/2})\biggr),
  \label{f^n}
\end{align}
where
\[
|\theta_2(un^{-1/2})|\le |\theta_1(un^{-1/2})| +\frac{e |u|^3}{2h_1^{ 6}n^{1/2}}.
\]

Furthermore, since $\psi (\lambda) = 1 - i/(2\lambda) + O(\lambda^2)$ as $\lambda\to 0,$ one has
\begin{align*}
\widehat{\psi} (un^{-1/2})
& =  \frac{1-e^{-iu \de n^{-1/2}} }{iu \de n^{-1/2}}
  \cdot  \frac{1-e^{-iu \D n^{-1/2}}}{iu \D n^{-1/2}}
 \notag
 \\
 & = 1 - \frac{iu (\D + \de)}{2n^{1/2}}
   + O\bigl( u^2  \D^2  n^{-1}\bigr),
\end{align*}
provided that $\de \le \D$.
Substituting the obtained  representations for $\f^n$ and $\widehat{\psi} $ into~\eqref{I_1} yields
\begin{align}
I_1  =
 n^{-1/2}  \int_{A_n} e^{-iuv - u^2/2} du
 & - \frac{i\mu_3}{6n} \int_{A_n}  u^3 e^{-iuv - u^2/2} du
 \notag\\
  &  - \frac{i(\D + \de)}{2n} \int_{A_n}  u  e^{-iuv - u^2/2} du
  +R^{(1)}_n,
  \label{I_11}
\end{align}
where it is not hard to show that, for $\D\le \D_1$ with a fixed $\D_1 >0$, one has the uniform bound
\begin{equation}
|R^{(1)}_n|\le   c(r,b)  n^{-(r-1)/2};
\label{R'}
\end{equation}
here and in what follows, by $c(r,b)\in (0,\infty) $ we denote  constants that can depend on the values of $r$ and $b$, and may be different even within one and the same formula.   To indicate how the bound~\eqref{R'}   was obtained, just note that
\begin{align*}
\biggl|n^{-1/2}
  & \int_{A_n} e^{-iuv - u^2/2}  \frac{u^3}{ n^{1/2}}\theta_2(un^{-1/2})du
\biggr|
\le
n^{-1} \int_{A_n} e^{ - u^2/2} |u|^3 |\theta_2 (un^{-1/2})|du
 \\
&
\le
n^{-1}   \int_{A_n} e^{ - u^2/2} |u|^3
 \biggl(|\theta_1 (un^{-1/2})| + \frac{e|u|^3}{2h_1^6 n^{1/2}}\biggr)du
 \\
 &
 \le
 n^{-1}   c (r,b)   \int  e^{ - u^2/2} |u|^3 \biggl( \frac{|u|^{r-3}}{n^{(r-3)/2}} +
  \frac{ |u|^3}{  n^{1/2}}\biggr) du
 \le c(r,b) n^{-(r-1)/2}.
\end{align*}
The contributions of other cross-products  in the expression for $\f^n \widehat{\psi} $ are bounded in a similar way, replacing $|\mu_3|$ with its upper bound $b^{3/r}$ due to Lyapunov's moment  inequality.

Clearly, replacing in~\eqref{I_11} the integrals $\int_{A_n}$ with $\int_{\R}$ will introduce an error of the order $o(n^{-2})$ uniform over the class ${\mathcal F}_{r,b}$. Therefore, making that change and then   replacing  the integrals over the whole line with their explicit values, and setting $\de:=\D n^{-1}$, we conclude that
\begin{equation}
I_1 = 2\pi   \phi (v)\biggl(n^{-1/2}
+ \frac{ \mu_3}{6n} (v^3 - 3v) - \frac{\D  }{2n} v\biggr) + R^{(2)}_n,
\label{I_1b}
\end{equation}
where for $R^{(2)}_n$ holds true the same bound~\eqref{R'} as for $R^{(1)}_n$.

To bound $I_2,$  note that it follows from~\eqref{theta} and \eqref{la^2} that
\[
 | (1-\f (\lambda)) -\lambda^2 /2|\le  \frac12 b^{3/r-1} {\lambda^2} \quad\mbox{for}\quad
  |\lambda|< h_2: =b^{-1}.
 \]
Therefore, setting $g:= \frac12 (1-b^{3/r-1})\in (0,\frac12),$ we have
\begin{align*}
|\f (\lambda) |
 & =|1-  (1-\f (\lambda) -\lambda^2 /2)- \lambda^2 /2|
\notag \\
 &
 \le |1-g  \lambda^2| \le e^{-g\lambda^2 }, \qquad  |\lambda|< h_2,
\end{align*}
implying that
\begin{align}
|I_2|
  & \le  \int_{B_2} | e^{-i\lambda x} \f^n (\lambda) \widehat{\psi} ( \lambda) | d\lambda
  \le
  \int_{B_2} |  \f (\lambda)   |^n d\lambda
  \notag
   \le
  \int_{|\lambda|\ge h_1 n^{-1/3}}  e^{-n g\lambda^2 }  d\lambda
    \\
  &= n^{-1/2}\int_{|u|\ge h_1 n^{1/6}}  e^{-gu^2 }  du
   \le
   n^{-1/2} \biggl (1 + \frac{1}{2 g h_1^2 n^{1/3}}\biggr)e^{- g h_1^2 n^{1/3} }.
   \label{I_2}
\end{align}

It remains to bound $I_3$. Since $|\widehat{\psi} (\lambda)|\le 2^2/(\Delta\delta\lambda^2) ,$ $\delta= \Delta n^{-1}$ and $h_2=b^{-1}$, one has
\begin{align*}
|I_3|
  & \le \rho^n
  \int_{|\lambda|\ge h_2} |\widehat{\psi} (\lambda)|  d\lambda
   \le \frac{4  \rho^n n }{\Delta^2}
   \int_{|\lambda|\ge 1/b}\frac{ d\lambda}{\lambda^2}
   = 8 b  \Delta^{-2} \rho^n n  .
\end{align*}
Choosing an arbitrary fixed  $q\in (\rho^{1/2},1)$ and setting $\eta:=\rho q^{-2}\in (0,1)$, we  will have
\begin{align}
|I_3| \le  8 b  \eta^n n  \quad \mbox{for all}\quad \Delta \ge q^n.
\label{I_3}
\end{align}

Now returning to~\eqref{iii} and using representation~\eqref{I_1b} for $I_1$ and the bounds     \eqref{I_2},  \eqref{I_3} for  $I_2$ and $I_3$,  we obtain that, for $\Delta\in [q^n, \Delta_1]$ with an arbitrary fixed $\Delta_1  \ge 0,$ one has
\begin{align}
\pr \bigl(\widetilde{S}_n \in \D [x)\bigr)
 & =  \D \phi (v)
  \biggl(n^{-1/2}
+ \frac{ \mu_3}{6n} (v^3 - 3v) - \frac{\D  }{2n} v\biggr) + \D R^{(3)}_n,
\label{for_smoothed}
\end{align}
where for $R^{(3)}_n$ holds  true the same bound~\eqref{R'}  as for  $R^{(1)}_n$.

Now set $ (\D-\delta) [x):= [x, x+\Delta- \delta)$. Clearly, $\{\widetilde{S}_n \in (\D-\delta) [x)\}\subset \{ {S}_n \in \D [x)\}$, so that  from~\eqref{for_smoothed} we see that
\begin{align}
\pr \bigl( S_n & \in \D  [x)\bigr)
   \ge \pr \bigl(\widetilde{S}_n \in (\D-\delta)  [x)\bigr)
 \notag \\
 &=
 \D(1-n^{-1}) \phi (v)  \biggl(n^{-1/2}
+ \frac{ \mu_3}{6n} (v^3 - 3v) - \frac{\D(1-n^{-1})  }{2n} v\biggr)
 + \D R^{(3)}_n
 \notag
 \\
 & =   \D \phi(v) \biggl(n^{-1/2}
+ \frac{ \mu_3}{6n} (v^3 - 3v) - \frac{\D   }{2n} v\biggr) + \D R^{(4)}_n,
 \label{lo_bound}
\end{align}
where for $R^{(4)}_n$ holds  true the same upper bound~\eqref{R'}  as for  $R^{(1)}_n$.

Setting $\overline{S}_n:= S_n + \delta U$, repeating the above calculations and using the obvious inclusion $ \{ {S}_n \in \D [x)\}\subset \{\overline{S}_n \in (\D+\delta) [x)\},$ we obtain an upper bound for the probability $\pr \bigl( S_n   \in \D  [x)\bigr)$ which is of the same form as the right-hand side of~\eqref{lo_bound}. Since clearly $n^{(r-3)/2}\cdot nR_n^{(4)}<c (r,b)<\infty$ according to the bound~\eqref{R'}, that completes the proof of part~(ii) of our theorem.

 \medskip

{\bf Acknowledgments.} The authors are grateful to the anonymous referees for drawing their attention to the note~\cite{Fo82}.


\begin{thebibliography}{X}

\bibitem{BoMo99}
A.~A.~Borovkov and A.~A.~Mogulskii. (1999)
Integro-local limit theorems includeing large deviationsfor for sums of random vectors.~I.
{\em Theory Probab.\ Appl.} 43, 1--12

\bibitem{BoMo01}
A.~A.~Borovkov and A.~A.~Mogulskii. (2001)
Integro-local limit theorems includeing large deviationsfor for sums of random vectors.~II.
{\em Theory Probab.\ Appl.} 45, 3--22.

\bibitem{Bo11}
A.~A.~Borovkov. (2011)
Integro-Local and Local Theorems on Normal and Large Deviations of the Sums of Nonidentically Distributed Random Variables in the Triangular Array Scheme. {\em Theory Probab.\ Appl.} 54, 571--587.

\bibitem{Bo13}
A.~A.~Borovkov. (2013)
{\em Probability theory.}
Springer, London.


\bibitem{Bo16}
A.~A.~Borovkov. (2016)
On the distribution of the first passage time of an arbitrary remote boundary by a random walk. {\em Theory Probab.\ Appl.} 61, 1--24.


\bibitem{Bo16a}
A.~A.~Borovkov. (2016)
Generalization and refinement of the integro-local Stone theorem for sums of random vectors. Submitted to  {\em Theory Probab.\ Appl.}






\bibitem{Ch87}
P.~L.~Chebyshev. (1890)
Sur deux th\'eor\`emes relatifs aux probabilit\'es.
{\em Acta Math.}  14, 305--315.




\bibitem{Cr28}
H.~Cram\'er. (1928)
On the composition of elementary errors. First paper: Mathematical deductions. {\em
Scandinavian Actuarial Journal},   1, 13--74.

\bibitem{Ed05}
F.~Y.~Edgeworth. (1905)
The Law of Error.
{\em Cambridge Philos.\ Soc.} 20, 36--66 and 113--141.


\bibitem{Es45}
C.~G.~Esseen. (1945)
Fourier analysis of distribution functions. A mathematical study of the Laplace-Gaussian law. {\em Acta Math.} 77, 1--125.

\bibitem{Fe67}
W.~Feller. (1971)
{\em An Introduction to Probability Theory and Its Applications.} Vol.~2, 2nd edn. Wiley, NY.


\bibitem{Fo82}
S.~V.~Fomin (1982)
Asymptotic expansions in the multidimensional central limit theorem. (Russian. English summary).
{\em Vestnik Leningrad.\ Univ.\ Mat.\ Mekh. Astronom.} no.~2, 116--118, 128.

\bibitem{GnKo49}
B.~V.~Gnedenko and A.~N.~Kolmogorov. (1954)
{\em Limit distributions for sums of independent random variables.} Addison-Wesley, Cambridge, Mass. [Translated from the 1949 Russian edition.]


\bibitem{Lo55}
M.~Loeve. (1955).
{\em Probability Theory.} D.~Van Nostrand Company, NY.



\bibitem{Pe95}
V.~V.~Petrov. (1995)
{\em Limit Theorems of Probability Theory: Sequences of Independent Random Variables.} Clarendon Press, Oxford.  

\bibitem{Rv54}
E.~L.~Rva\v ceva. (1962)
On domains of attraction of multi-dimensional distributions.
In: {\em Select.\ Transl.\ Math.\ Statist.\ and Probability,} Vol.~2, American Mathematical Society, Providence, R.I., 183--205. [Trnaslated from the original  1954 Russian paper.]

\bibitem{Sh64}
L.~A.~Shepp (1964)
A local limit theorem.
{\em Ann.\ Math.\ Statist.} 35,  419--423.

\bibitem{St65}
C.~Stone. (1965)
A local limit theorem for nonlattice multi-dimensional distribution functions.
{\em Ann.\ Math.\ Statist.} 36, 546--551.

\bibitem{St67}
C.~Stone. (1967)
On local and ratio limit theorems.
In: {\em  Proc.\ Fifth Berkeley Sympos.\ Math.\ Statist.\ and Probability } (Berkeley, Calif., 1965/66), Vol. II: Contributions to Probability Theory, Part 2, pp.~217--224. Univ.\ California Press, Berkeley, CA.


\end{thebibliography}
\end{document}